\documentclass[12pt]{amsart}

\usepackage{amsmath}
\usepackage{amssymb}
\usepackage{amsfonts}
\usepackage{amsthm}
\usepackage{enumerate}
\usepackage{hyperref}
\usepackage{color}
\usepackage{psfrag}
\usepackage[all]{xy}
\usepackage{comment}

\theoremstyle{plain}
\newtheorem{thm}{Theorem}

\newtheorem{prop}[thm]{Proposition}
\newtheorem{remark}[thm]{Remark}

\theoremstyle{definition}

\newtheorem{example}[thm]{Example}

\def\projdim{\mathop{\rm projdim}\nolimits}
\def\reg{\mathrm{reg}}
\def\reg{\mathop{\rm reg}\nolimits}

\begin{document}


\title[Splitting and regularity of powers of monomial ideals]{Splitting and regularity of powers of monomial ideals}
\title[Differences in regularities]{Differences in regularities of a power
and its integral closure and symbolic power}


\author[]{Irena Swanson}

\author[]{Siamak Yassemi}

\address{Department of Mathematics, Purdue University, West Lafayette,
IN 47907}

\email{irena@purdue.edu, syassemi@purdue.edu}


\begin{abstract}
We show that for a square-free monomial ideal,
the regularity of its symbolic (second) power and of the integral closure of
of its (second) power
can differ from the regularity of its ordinary (second) power
by an arbitrarily large integer.

\end{abstract}






\thanks{}


\maketitle


All ideals in this paper are homogeneous in polynomial rings over fields,
and for any ideal~$I$ and any integer~$s$,
$\overline{I^s}$ denotes the integral closure of~$I^s$
and $I^{(s)}$ denotes the $s^{th}$ symbolic power of~$I$.
For monomial ideals we also use the $s$th (special)
Frobenius power $I^{[s]}$ of~$I$
to be the ideal generated by the $s$th powers of all the monomials in~$I$.
For any square-free monomial ideal~$I$ and any integer $s\ge 1$,
$I^{[s]} \subseteq I^s \subseteq \overline{I^s} \subseteq I^{(s)}$.

Kumar and Kumar~\cite{KK} showed on a monomial example~$I$
that for all positive integers~$s$,
the regularity of $I^s$ can differ from
the regularity of $\overline{I^s}$ and the regularity of~$I^{(s)}$ by~$1$,
and they showed on a square-free monomial example~$I$
that the regularity of $I^2$ can differ from
the regularity of $\overline{I^2}$ and
the regularity of $I^{(2)}$ by~$1$.
Most prior work in this area has been on establishing
asymptotic properties and upper and lower bounds
on the regularities of~$I^s$, $\overline{I^s}$ and $I^{(s)}$,
as well as on finding cases of ideals for which
the regularities of these ideals agree.
See Kumar and Kumar~\cite{KK}
and Minh, Nam, Phong, Thuy and Vu~\cite{MNPTV}
for a summary of known results of cases of edge ideals
for which the regularity of a power
equals the regularity of the integral closure of the same power,
and for the very few known classes of edge ideals
for which the regularity of a power
equals the regularity of that symbolic power.


For the rest of the paper
we fix the square-free ideal of Sturmfels from~\cite{S}:
$$
J = (x_1x_4x_5, x_1x_3x_6, x_2x_3x_4, x_2x_5x_6, x_3x_4x_5,
x_3x_4x_6, x_3x_5x_6, x_4x_5x_6)
$$
in the polynomial ring $S = k[x_1, \ldots, x_6]$ in variables
$x_1, \ldots, x_6$ over a field~$k$.
Sturmfels showed that
$J$ is an 8-generated square-free monomial ideal with a linear resolution
independent of the characteristic,
and that $J^2$ does not have a linear resolution.
Kumar and Kumar showed in~\cite{KK} that
$\reg(J^2) = 7$ and $\reg(\overline{J^2}) = \reg(J^{(2)}) = 6$.
Using the graded Betti numbers of $J^{[2]}, J^2, \overline{J^2}$ and $J^{(2)}$
and using the splitting homomorphisms from~\cite{KS},
we construct in this paper square-free monomial ideals~$I$
in polynomial rings
for which the differences
$$
\reg(I^2) - \reg(\overline{I^2}),
\reg(\overline{I^2}) - \reg(I^{[2]}),
\reg(I^{[2]}) - \reg(I^{(2)})
$$
can be arbitrarily large.
More specifically,
we prove the following:

\begin{thm}
\label{thmmain}%
Let $e$ be a positive integer,
let $R$ be the polynomial ring
$k[y_{i,j} : i = 1, \ldots, 6; j = 1, \ldots, e]$,
and let $\varphi_e : S \to R$ be the $k$-algebra homomorphism
that for each $i$
takes $x_i$ into the product $y_{i,1} y_{i,2} \cdots y_{i,e}$.
Set $I_e = \varphi_e(J) S$.
Then
\begin{enumerate}
\item
$I_e$ is a square-free monomial ideal.
\item
$I_e^{[2]} = \varphi_e(J^{[2]})$,
$\reg(I_e^{[2]}) = 10e -2$.
\item
$I_e^2 = \varphi_e(J^2)$,
$\reg(I_e^2) = 12e-5$.
\item
$\overline{I_e^2} = \varphi_e(\overline{J^2})$,
$\reg(\overline{I_e^2}) = 11e - 5$.
\item
$I_e^{(2)} = \varphi_e(J^{(2)})$,
$\reg(I_e^{(2)}) = 9e - 3$.
\end{enumerate}
\end{thm}

\begin{proof}
It is very easy to compute the generators of $J^{[2]}$ and $J^2$.
It is straightforward to compute $J^{(2)}$
by first establishing the minimal primes
$(x_3, x_5)$,
$(x_4, x_6)$,
$(x_1, x_3, x_6)$,
$(x_1, x_4, x_5)$,
$(x_2, x_3, x_4)$,
$(x_2, x_5, x_6)$
over $J$,
then establishing the corresponding primary components of $J^2$,
and finally taking the intersection of the six primary ideal.
We get that
\begin{align*}
J^{(2)}
= &(x_3x_4x_5x_6, x_4^2x_5^2x_6^2, x_2x_4x_5^2x_6^2, x_3^2x_5^2x_6^2,
x_2x_3 x_5^2x_6^2, x_2^2x_5^2x_6^2, x_1x_3^2x_5x_6^2, x_1x_2x_3x_5x_6^2,\\
& x_3^2x_4^2x_6^2, x_1x_3^2x_4x_6^2,
x_1^2x_3^2x_6^2, x_1x_4^2x_5^2x_6, x_1x_2x_4x_5^2x_6,
x_2x_3^2x_4^2x_6, x_1x_2x_3^2x_4x_6,\\
& x_3^2x_4^2x_5^2, x_1x_3x_4^2x_5^2, x_1^2x_4^2x_5^2, x_2x_3^2x_4^2x_5,
x_1x_2x_3x_4^2x_5, x_2^2x_3^2x_4^2).
\end{align*}
According to Macaulay2,
the integral closure of $J^2$ is
\begin{align*}
\overline{J^2}
=
&(x_4^2x_5^2x_6^2, x_3x_4x_5^2x_6^2, x_2x_4x_5^2x_6^2, x_3^2x_5^2x_6^2,
x_2x_3x_5^2x_6^2, x_2^2x_5^2x_6^2, x_3x_4^2x_5x_6^2, x_3^2x_4x_5x_6^2,\\
& x_2x_3x_4x_5x_6^2, x_1x_3x_4x_5x_6^2, x_1x_3^2x_5x_6^2, x_1x_2x_3x_5x_6^2,
x_3^2x_4^2x_6^2, x_1x_3^2x_4x_6^2, x_1^2x_3^2x_6^2, \\
& x_3x_4^2x_5^2x_6,
x_1x_4^2x_5^2x_6, x_3^2x_4x_5^2x_6, x_2x_3x_4x_5^2x_6, x_1x_3x_4x_5^2x_6,
x_1x_2x_4x_5^2x_6, x_3^2x_4^2x_5x_6,\\
& x_2x_3x_4^2x_5x_6, x_1x_3x_4^2x_5x_6,
x_2x_3^2x_4x_5x_6, x_1x_3^2x_4x_5x_6, x_2^2x_3x_4x_5x_6,
x_1x_2x_3x_4x_5x_6,\\
& x_1^2x_3x_4x_5x_6, x_2x_3^2x_4^2x_6, x_1x_2x_3^2x_4x_6, x_3^2x_4^2x_5^2,
x_1x_3x_4^2x_5^2, x_1^2x_4^2x_5^2, x_2x_3^2x_4^2x_5,\\
& x_1x_2x_3x_4^2x_5, x_2^2x_3^2x_4^2).
\end{align*}
It is not hard to verify that all the listed monomials are
integral over $J^2$;
for the remainder of the paper we assume that the integral closure
is as stated.

For monomial ideals it is doable
to establish free resolutions by hand.
The Betti tables below were provided by Macaulay2
(we do not write out the all-zero rows).
The Betti table for $J^{[2]}$ is:
$$
\begin{matrix}
& 0 & 1 & 2 \cr
\hbox{6:} & 8 & \cdot & \cdot \cr
\hbox{7:} & \cdot & 11 & \cdot \cr
\hbox{8:} & \cdot & \cdot & 4 \cr
\end{matrix}
$$
for $J^2$ it is:
$$
\begin{matrix}
& 0 & 1 & 2 & 3 & 4 & 5 \cr
\hbox{6:} & 36 & 84 & 75 & 32 & 6 & \cdot \cr
\hbox{7:} & \cdot & 1 & 4 & 6 & 4 & 1 \cr
\end{matrix}
$$
for $\overline{J^2}$ it is:
$$
\begin{matrix}
& 0 & 1 & 2 & 3 & 4 & 5 \cr
\hbox{6:} & 37 & 90 & 89 & 48 & 15 & 2 \cr
\end{matrix}
$$
and for $J^{(2)}$ it is:
$$
\begin{matrix}
& 0 & 1 & 2 & 3 \cr
\hbox{4:} & 1 & \cdot & \cdot & \cdot \cr
\hbox{6:} & 20 & 40 & 24 & 4 \cr
\end{matrix}
$$
For readers not familiar with Macaulay2's notation of Betti tables
we read off the resolution for the last ideal:
$$
0 \to R(-6-3)^4
\rightarrow R(-6-2)^{24}
\rightarrow R(-6-1)^{40}
\rightarrow R(-4-0) \oplus R(-6-0)^{20}
\rightarrow J^{(2)} \to 0.
$$

We now apply $\varphi_e$.
Certainly $I_e = \varphi_e(J)$ is square-free and monomial,
proving~(1).

Say by \cite[Theorem 1.2]{KS},
$\varphi_e$ is a flat homomorphism,
and so it preserves all Betti numbers.
By the homogeneous nature of $\varphi_e$,
the degrees in the resolution of $\varphi_e(J)$
are $e$ times the degrees in the resolution of $J$.
Thus we can read off the following:
\begin{align*}
&\reg(\varphi_e(J^{[2]})) = \max\{6e, 8e-1, 10e-2\} = 10e -2, \\
&\reg(\varphi_e(J^2)) = \max\{6e, 8e-1, 9e-2, 10e-3, 11e-4, 12e-5\} = 12e-5, \\
&\reg(\varphi_e(\overline{J^2})) = \max\{6e, 7e-1, 8e-2, 9e-3, 10e-4, 11e-5\}
= 11e - 5, \\
&\reg(\varphi_e(J^{(2)})) = \max\{6e, 7e-1, 8e-2, 9e-3\} = 9e -3.
\end{align*}

It remains to prove the first equalities in parts (2) through~(4).
The equalities in~(2) and~(3) hold by the definition of homomorphisms.
The equality in~(4) holds by Proposition~\ref{propphicommuteswithpowers} below.
It was proved in~\cite[Lemmas 1.5, 1.6]{KS}
that $\varphi_e$ takes irredundant primary decompositions
to irredundant primary decompositions,
associated primes to associated primes,
while preserving heights and 
inclusions and non-inclusions of associated primes.
Thus the equality holds in~(5).
\end{proof}

\begin{prop}
\label{propphicommuteswithpowers}%
Let $J$ be a monomial ideal in a polynomial ring~$S$ over a field~$k$.
Let~$R$ be another polynomial ring over~$k$
and let $\varphi : S \to R$ be a $k$-algebra homomorphism
that takes each of the variables of $S$
into a product of distinct and disjoint variables in~$R$.
More precisely,
$\varphi(x_i) = v_{i,1} v_{i,2} \cdots v_{i,t_i}$
where $x_1, \ldots, x_u$ are the variables in $S$,
$t_1, \ldots, t_u$ are positive integers,
and the $t_1 + t_2 + \cdots + t_u$ variables $v_{i,j}$ in~$R$ are distinct.
Then for any positive integer~$s$,
$$
\overline{\varphi(J^s) R} = \varphi(\overline {J^s}) R.
$$
\end{prop}

\begin{proof}
For any ring homomorphism $\psi$ and any ideal~$L$ in the domain,
it is true that $\psi(\overline L) \subseteq \overline{\psi(L) R}$.
In particular,
$\varphi(\overline {J^s}) \subseteq \overline{\varphi({J^s}) R}$.

The integral closure of $\overline{\varphi(J^s) R}$ is generated by monomials~$m$
with the property that for some integer~$r$,
$m^r \in \varphi(m_1) \cdots \varphi(m_r) R$
$= \varphi(m_1 \cdots m_r) R$
for some generators $m_1, \ldots, m_r$ of~$J^s$ (page 9 in \cite{HS}).
By the nature of $\varphi$,
for each $i = 1, \ldots, n$,
the exponent of $v_{i,j}$
in $\varphi(m_1 \cdots m_r)$ is independent of $j = 1, \ldots, t_i$,
and so by the nature of the inclusion
$m^r \in \varphi(m_1 \cdots m_r) R$,
by possibly lowering some exponents in~$m$,
we may assume that the exponents
of $v_{i,1}, v_{i,2}, \ldots, v_{i,t_i}$ are all the same also in~$m$.
Thus $m$ can be written as $\varphi(n)$ for some $n \in R$.
Then we have $\varphi(n^r) = (\varphi(n))^r = m^r \in \varphi(m_1 \cdots m_r) R$.
By \cite[Theorem 1.2]{KS},
$\varphi$ is faithfully flat,
so that $n^r \in m_1 \cdots m_r S$.
This proves that $n \in \overline {J^s}$,
and so $m = \varphi(n) \in \varphi(\overline {J^s})R$.

This finishes the proof that
$\varphi(\overline {J^s})R = \overline{\varphi({J^s}) R}$.
\end{proof}

The homomorphisms in Proposition \ref{propphicommuteswithpowers}
were studied in~\cite{KS} and were called {\bf splittings}.
The homomorphisms $\varphi_e$ in Theorem~\ref{thmmain}
are very special splittings.
As pointed out in the proof of Theorem~\ref{thmmain},
if $\{\beta_{i,j}\}$ are the graded Betti numbers
of an (arbitrary) ideal~$I$ in~$S$,
with $i = 0, \ldots, \projdim(I)$,
then the graded Betti numbers of $\varphi_e(I)$ are $\{e\beta_{i,j}\}$.
In particular, if $d(I) = \max\{\beta_{0,j}: j\}$ denotes the maximal degree
of an element in a minimal generating set,
we have that $d(\varphi_e(I)) = e d(I)$,
and even more in particular,
$d((\varphi_e(I))^s) = e  d(I^s)$,
$d(\overline{(\varphi_e(I))^s}) = e d(\overline{I^s})$,
$d(\overline{(\varphi_e(I))^{(s)}}) = e d(\overline{I^{(s)}})$,
and $d(\overline{(\varphi_e(I))^{[s]}}) = e d(\overline{I^{[s]}})$.


Macaulay2 computed the regularities in the table below
for the splittings of Sturmfels's ideal~$J$,
when we split the variables listed in the first column into a product of $m$
distinct variables each
and the rest of the variables remain the product of one variable,
where $m$ is the entry in the second column.
Note that always $\varphi(J)^{[2]}
\subseteq \varphi(J)^2
\subseteq \overline{\varphi(J)^2}
\subseteq \varphi(J)^{(2)}$,
but that the regularity of the first ideal
may be strictly larger or strictly smaller than
the regularity of the second ideal on this list.
\def\vr{\smash{\vrule height 2.1ex depth .9ex width .01em}}
$$
\displayindent=1mm
\halign{#\hskip0.2em\vr&& \hskip0.2em \hfil # \hfil \hskip0.2em \cr
\hbox{variables to split\hfil} & into \#vars
& $\reg(\varphi(J)^{[2]})$
& $\reg(\varphi(J)^2)$
& $\reg(\overline{\varphi(J)^2})$
& $\reg(\varphi(J)^{(2)})$
\cr
\noalign{\vskip2pt\hrule\vskip2pt}
all & 1 each & 8 & 7 & 6 & 6 \cr
$x_1$ & 2 & 10 & 9 & 8 & 8 \cr
$x_1, x_2$ & 2 & 10 & 11 & 9 & 8 \cr
$x_1, x_2, x_3$ & 2 & 12 & 13 & 11 & 10 \cr
$x_1, x_2, x_3, x_4$ & 2 & 14 & 15 & 13 & 12 \cr
$x_1, x_2, x_3, x_4, x_5$ & 2 & 16 & 17 & 15 & 14 \cr
$x_2, x_3, x_4, x_5, x_6$ & 3 & 28 & 27 & 26 & 22 \cr
$x_3, x_4, x_5, x_6$ & 3 & 24 & 23 & 22 & 20 \cr
$x_4, x_5, x_6$ & 3 & 20 & 19 & 18 & 18 \cr
$x_5, x_6$ & 3 & 16 & 15 & 14 & 14 \cr
$x_6$ & 3 & 12 & 11 & 10 & 10 \cr
$x_1$ & 3 & 12 & 11 & 10 & 10 \cr
$x_1, x_2$ & 3 & 12 & 15 & 12 & 10 \cr
$x_1, x_2, x_3$ & 3 & 16 & 19 & 16 & 14 \cr
$x_1, x_2, x_3, x_4$ & 3 & 20 & 23 & 20 & 18 \cr
$x_1, x_2, x_3, x_4, x_5$ & 3 & 24 & 27 & 24 & 22 \cr
$x_1, x_3, x_5$ & 3 & 20 & 19 & 18 & 16 \cr
$x_2, x_4, x_6$ & 3 & 20 & 19 & 18 & 16 \cr
$x_1, x_3$ & 3 & 16 & 15 & 14 & 14 \cr
$x_4, x_6$ & 3 & 16 & 15 & 14 & 14 \cr
}
$$




\end{document}